\documentclass[12pt]{amsart}

\usepackage{amsmath,graphics,comment}
\usepackage{amsfonts,amssymb}

  \includecomment{inlong} \excludecomment{inshort}   
 \includecomment{inlow} \excludecomment{inhigh}

\theoremstyle{plain}
\newtheorem*{theorem*}{Theorem}
\newtheorem*{lemma*} {Lemma}
\newtheorem*{corollary*} {Corollary}
\newtheorem*{proposition*} {Proposition}
\newtheorem{theorem}{Theorem}[section]
\newtheorem{lemma}[theorem]{Lemma}
\newtheorem{corollary}[theorem]{Corollary}
\newtheorem{proposition}[theorem]{Proposition}
\newtheorem{conjecture}[theorem]{Conjecture}

\theoremstyle{remark}
\newtheorem*{remark}{Remark}

\newtheorem*{claim}{Claim}

\theoremstyle{definition}

\def \R {\mathbf{R}}
\def \Z {\mathbf{Z}}
\def \C {\mathbf{C}}

\def\eps{\epsilon}

\def\s{\sigma}

\def\id{\mbox{id}}

\def\sign{\mbox{sign}}
\def\Z{\Bbb{Z}}
\def\C{\Bbb{C}}

\def\N{\Bbb{N}}
\def\l{\lambda}
\def\part{\partial}

\def\a{\alpha}

\def\bp{\begin{pmatrix}}

\def\sm{\setminus}
\def\ep{\end{pmatrix}}
\def\bn{\begin{enumerate}}

\def\en{\end{enumerate}}
\def\ba{\begin{array}}
\def\ea{\end{array}}

\def\S{\Sigma}
\def\s{\sigma}
\def\a{\alpha}
\def\b{\beta}
\def\w{\omega}
\def\ti{\tilde}
\def\lk{\mbox{lk}}
\def\vph{\varphi}
\def\fr12{\frac{1}{2}}
\def\diag{\mbox{diag}}

\def\Aut{\mbox{Aut}}

\def\im{\mbox{Im}}
\def\Im{\mbox{Im}}

\def\t{\theta}

\def\v{\varphi}

\def\mapstoleft{\leftarrow \hspace{-0.12cm} \line(0,1){5.3}}
\def\cmtbf#1{} \def\cmt#1{}

\begin{document}

\title{Link concordance, boundary link concordance and eta-invariants}
\author{Stefan Friedl}
\date{\today}
\begin{abstract}
We study the eta-invariants of links and show that in many cases they form link concordance invariants, in 
particular that many eta-invariants vanish for slice links.
This result contains and generalizes previous invariants by Smolinsky and Cha--Ko.
We give a formula for the eta-invariant for boundary links. In several intersting cases this allows us 
to show that a given link is not slice. We show that even
more eta-invariants have to vanish for boundary slice links.  
We give an example of a boundary link $L$ that is not boundary slice but where all the known link concordance invariants
computed so far are zero.
\end{abstract}
\maketitle

\section{Introduction}
An $m$-link of dimension $n$ is an embedded oriented smooth submanifold of $S^{n+2}$ that is homeomorphic to
 $m$ ordered copies of $S^{n}$. A link concordance
between two given links in $S^{n+2}$ is a properly embedded oriented submanifold in $S^{n+2} \times [0,1]$ that is
homeomorpic to $m$ copies of $S^n \times [0,1]$ and intersects $S^{n+2} \times 0$ and  $S^{n+2} \times 1$ at the given links. 
We say a link is slice if it is concordant to the trivial link. Equivalently a link is  slice if it bounds
$m$ disjoint smooth disks in $D^{n+3}$.

Denote by $C(n,m)$ the set of concordance classes of $m$-links of dimension $n$.
The set $C(n,1)$ is just the set of knot concordance classes,
it has a well-defined group structure given by
connected sum along arcs. 
Connected sum of links does not give a well-defined group structure on $C(n,m)$ since
there's no canonical choice of arcs (cf. proposition
\ref{propnogroup}).  

It is very difficult to determine $C(n,m)$, a common approach is
to study links with some extra structure. A boundary link is an
$m$-link which has $m$ disjoint Seifert manifolds, i.e. there
exist $m$ disjoint oriented $(n+1)$-submanifolds $V_1,\dots,V_m
\subset S^{n+2}$ such that $\partial(V_i)=L_i, i=1,\dots,m$. A
boundary link concordance between two given boundary links in
$S^{n+2}$ is a link concordance which  bounds $m$ disjoint
$(n+2)$-manifolds in $S^{n+2}\times [0,1]$. We say $L$ is boundary
slice if it is boundary concordant to the unlink. Denote by
$B(n,m)$ the set of boundary concordance classes of $m$-boundary
links of dimension $n$.

A pair $(L,V)$ consisting of a boundary link and a Seifert manifold is called boundary link pair.
There's an obvious notion of concordance for boundary link pairs.
Denote by $C_n(B_m)$ the set of concordance classes of boundary link pairs.
Assume $n>1$ and let $(L_1,V_1)$ and $(L_2,V_2)$ be representatives of elements in $C_n(B_m)$.
We can assume that $V_1,V_2$ are simply connected, and then $S^{n+2}\sm (V_1\cup V_2)$ is simply connected,
in particular up to homotopy there's only one arc connecting $L_1,L_2$ in the complement
of $V_1 \cup V_2$, therefore the boundary connected sum $(L_1\# L_2,V_1 \# V_2)\in C_n(B_m)$ 
is well-defined and turns $C_n(B_m)$ into a group.

Let $F_m$
be the free group on the generators $t_1,\dots,t_m$. An $F_m$-link is a pair $(L,\v)$ 
where $L$ is a link in $S^{n+2}$ and
$\v:\pi_1(S^{2n+2}\sm L) \to F_m$ is an epimorphism sending an $i^{th}$ meridian to $t_i$.
A pair
$(N,\Phi)$ is an $F_m$-concordance between $(L_0,\v_0)$ and $(L_1,\v_1)$ if $M$ is a link concordance
between the links $L_0$ and $L_1$ and $\Phi:\pi_1(S^{n+2}\times [0,1] \setminus N) \to F_m$ is a map extending
$\v_0$ and $\v_1$ up to inner automorphisms (cf. \cite{CS80}).  
Denote  by $C_n(F_m)$ the set of $F_m$-concordance classes of $F_m$-links.
If $n>1$ then any element in $C_n(F_m)$ has a representative $(L,\varphi)$ such that
$\varphi$ is an isomorphism, this $\varphi$ defines canonical meridians for $L$, which can be used
to defined a well-defined group structure on $C_n(F_m)$.
   
By the transversality argument, such an epimorphism $\v$ gives a Seifert surface $V_{\v}$. 
Conversely, the existence of a
Seifert surface $V$ for $L$ produces such an epimorphism $\v_{V}$ by the Thom-Pontryagin construction. We'll freely go back
and forth between isotopy classes of boundary link pairs $(L,V)$ and $F_m$-links $(L,\v)$.
Similarly there's an equivalence between the respective concordances,
in particular $C_n(B_m) \cong C_n(F_m)$, which is a group isomorphism for $n>1$.

We say that $\v:\pi_1(S^{n+2}\sm L) \to F_m$ is a splitting map if it sends meridians to generators.
There's  in general
not a unique splitting map.
Denote by $CA_m$ the group of automorphisms of $F_m$ which send $t_i$ to a conjugate of $t_i$ for
each $i=1,\dots,m$.  
\begin{lemma} \label{lemmacam}
If $\v:\pi_1(S^{n+2}\sm L) \to F_m$ is a splitting map, then for any $\phi \in CA_m$ the map $\phi \circ
\v$ is a splitting map as well, and in fact all splitting maps are of the form $\phi \circ \v$ for some
$\phi \in CA_m$.
\end{lemma}

This means that we have an action of $CA_m$ on $C_n(F_m)$. The inner automorphisms of $F_m$ are
elements in $CA_m$ and act trivially on $C_n(F_m)$. We therefore define
$A_m$ to be the quotient group of
$CA_m$ by the inner automorphisms of $F_m$. We get  an action of $A_m$ on $C_n(F_m)$.
Denote by $\phi_{ij}:F_m \to F_m$ the map which sends $t_i$ to $t_jt_it_j^{-1}$ and $t_k$ to $t_k$ for
$k \ne i$. 
We quote the following proposition (cf. \cite{K84}, \cite{K87}).

\begin{proposition} \label{propa2trivial}
$CA_m$ (and in particular $A_m$) is generated by $\phi_{ij}$ for $i,j=1,\dots,m$ and $i \ne j$.
Furthermore the groups $A_1,A_2$ are trivial.  
\end{proposition}

Under the isomorphism $C_n(F_m)\cong C_n(B_m)$ the group $A_m$ also acts on $C_n(B_m)$, the action
of $\phi_{ij}$ on a Seifert surface has been described explicitely by Ko \cite{K87}.
Ko \cite{K87} furthermore showed that $A_m$ acts non-trivially on $C_n(B_m)$ and hence acts non-trivially on $C_n(F_m)$.

\begin{theorem}\cite{CS80}  
\[ B(n,m) \cong C_n(F_m)/A_m \cong C_n(B_m)/A_m\]
\end{theorem}

Cappell and Shaneson showed that $C_{2k}(F_m)=0$, i.e. all even dimensional boundary links are boundary slice.
 It is not known whether all even dimensional (boundary) links are slice.
We'll restrict ourselves from now on to odd-dimensional links. 

For $\eps=\pm 1$ we call $A=(A_{ij})_{i,j=1,\dots,m}$ an $\eps$-boundary link Seifert matrix of size $(g_1,\dots,g_m)$ if
$A$ is a matrix with entries
$A_{ij}$ which are $(2g_i \times 2g_j)$--matrices over $\Z$ such that
$A_{ij}=-\eps A_{ji}^t$ for $i \ne j$ and $\det(A_{ii}+\eps A_{ii}^t)=1$
(cf. \cite{L77}, \cite{K87}). 
We say that $A_{ij}$ is metabolic if there exists a block diagonal matrix $P=\diag(P_1,\dots,P_m)$ such
that each $P_iA_{ij}P_j^t$ is of the form
\[ \bp 0 & C \\ D&E \ep \]
where $0$ is a $g_i\times g_j$-matrix.
This generates in a natural way an equivalence class of matrices,
the set of equivalence classes is denoted by $G(m,\eps)$.

If $n=2q-1$ then picking a basis for the torsion free parts of $H_q(V)=H_q(V_1)\oplus \dots\oplus H_q(V_m)$ 
we can associate to a boundary link pair $(L,V)$ the matrix representing the Seifert pairing
\[ \ba{rcl} H_q(V)\times H_q(V) &\to &\Z \\
(a,b)&\mapsto &\lk(a,b_+)
\ea \]

\begin{theorem} \cite{K85} \label{theoremboundaryclass}
\bn
\item Every Seifert matrix is the Seifert matrix of a boundary link pair,
\item for $q \geq 3$ 
\[ C_{2q-1}(B_m) \cong G(m,(-1)^q) \]
\item $C_3(B_m)$ is isomorphic to a subgroup of $G(m,1)$ of index $2^m$.
\en
\end{theorem}
The $A_m$ action on  $C_{2q-1}(B_m)$ translates to an action of  $A_m$
on $G(m,(-1)^q)$ which was explicitely computed by Ko \cite{K87}.
Summarizing we get for $q\geq 3$ that
\[ B(2q-1,m)\cong C_{2q-1}(B_m)/A_m\cong G(m,(-1)^q)/A_m \]
 Levine \cite{L69b} showed that $G(1,\eps)\cong \Z^{\oplus\infty}\oplus
\Z_2^{\oplus\infty}\oplus\Z_4^{\oplus\infty}$ (cf. also \cite{S77}).
Recently Sheiham \cite{S02} showed  that for $m>1$,
 $G(m,\eps)\cong \Z^{\oplus\infty}\oplus
\Z_2^{\oplus\infty}\oplus\Z_4^{\oplus\infty}\oplus\Z_8^{\oplus \infty}$,
furthermore Sheiham defined full invariants for $G(m,\eps)$.


A lot of effort has been put into the study of the forgetful map
\[  
 B(n,m)\to C(n,m) \]
Cochran and Orr \cite{CO90}, \cite{CO93}, Gilmer and Livingston \cite{GL92} 
and Levine \cite{L94} showed that this map   is not
surjective, i.e. there exist links which are not concordant to boundary links.
It is an open question whether the kernel is trivial, i.e. whether any knot that is slice is also boundary
slice. It would be very difficult to find
counter-examples in dimension one, since one can easily see that any ribbon (boundary) link is boundary slice.
 

Given a closed smooth odd dimensional manifold $M$ and a unitary representation $\a:\pi_1(M)\to U(k)$,
Atiyah--Patodi--Singer
\cite{APS75} defined an invariant $\eta_{\a}(M)\in \R$, called the eta--invariant, which can be computed
in terms of signatures of bounding manifolds, if these exist.
 For a group $G$ a pair $(M,\v)$ is called a $G$-manifold if $M$ is a smooth odd-dimensional manifold  
and $\v:\pi_1(M)\to G$ a homomorphism. Define $\rho(M,\varphi):{R}_k(G) \to  \R$
via $\rho(M,\varphi)(\a):=\eta_{\a\circ \varphi}(M)$.
Two   $G$-manifolds $(M_j,\a_j), j=1,2$ are called
homology $G$-bordant if there exists a
$G$-manifold
$(N,\b)$  such that $\partial(N)=M_1 \cup -M_2, H_*(N,M_j)=0$ for $j=1,2$
and, up to inner automorphisms of $G$,
$\beta|\pi_1(M_j)=\a_j$.

\begin{theorem} \cite[p. 95]{L94}  
If  $(M_i,\a_i), i=1,2$ are homology $G$-bordant manifolds, then 
$\rho(M_1,\v_1)(\a)=\rho(M_2,\v_2)(\a)$  for all $\a:G\to U(k)$ that factor through a $p$--group.
  \end{theorem}

We'll study the $\rho$--invariant for $M_L$, the result of zero framed surgery along $L\subset S^{2q+1}$.
For $G$ a group define the lower central series inductively by $G_0:=G, G_i:=[G,G_{i-1}]$.
For the remainder of the introduction we'll denote the free group on $m$ generators by $F$.
For an $m$-component link $L\subset S^{2q+1}$ we have in many cases (e.g. if $q>1$) an isomorphism
$\pi_1( S^{2q+1}\sm L)/\pi_1( S^{2q+1}\sm L)_i \to F/F_i$. A choice of isomorphism is called an $F/F_i$--structure.
Two links $L_1,L_2$ with $F/F_i$--structures that are concordant also have 
concordant
$F/F_i$--structures, and 
 $M_{L_1}$ and $M_{L_2}$ have homology $F/F_i$--bordant $F/F_i$--structures. Applying the above theorem
gives a link concordance obstruction theorem. The theory becomes even easier if we want to find
sliceness obstructions since any slice knot has an $F/F_i$--structure for all $i$ and since any representation factoring
through a $p$--group factors through $F/F_i$ for some $i$.

\begin{theorem}  
Let $L\subset S^{2q+1}$ be a slice link, if $\a:\pi_1(M_L)\to U(k)$ factors through a $p$--group, then  
$\eta_{\a}(M_L)=0$.
  \end{theorem}

 Define $PD(k)\subset U(k)$ to be the subgroup generated by permutation matrices and diagonal matrices.
For a prime $p$ define $PD_p(k)\subset PD(k)$ to be the subgroup of matrices where all 
eigenvalues are roots of unity of order
a power of $p$.  
 
\begin{theorem} \label{bestslicethmintro}
Let $L\subset S^{2q+1}$ be a slice link with meridians $\mu_1,\dots,\mu_m$. Let $p$ be a prime number
and let $U_1,\dots,U_m\in PD_p(K)$. Then there exists a unique representation $\b:\pi_1(M_L)\to U(k)$ with
$\b(\mu_j)=U_j$. Furthermore $\eta_{\b}(M_L)=0$.
\end{theorem}

This  gives the best possible sliceness obstruction
theorem that can be based on Levine's theorem.
 These obstructions combine, simplify and generalize sliceness obstructions defined by Smolinsky
\cite{S89},
\cite{S89b} and Cha and Ko \cite{CK99}.

For an $F$--link $(L,\v)$ the $\rho$-invariant can be explicitely computed   in terms of its Seifert
matrix. In the case $n=4q+3$ the following holds, the case $n=4q+1$ being only marginally more complicated
(cf. theorem \ref{thmbdyetacomp}).

\begin{theorem} 
Let $(L\subset S^{4q+3},\v)$ be an $F_m$-link, 
$A=(A_{ij})_{i,j=1,\dots,m}$ a Seifert matrix,
 $\a:F_m\to U(k)$ a representation.
Let $U_i:=\a(t_i)$, then $\rho(M_L,\v)(\a)=\sign(M(A,\a))$ where $M(A,\a)$ equals
\[
 \bp A_{11}\otimes(\id-U_1^{-1})+A_{11}^t\otimes(\id-U_1 ) & A_{12}\otimes(\id-U_1)(\id-U_2^{-1}) &
\dots
\\ A_{21}\otimes(\id-U_2)(\id-U_1^{-1}) &A_{22}\otimes(\id -U_2^{-1})+ A_{22}^t\otimes(\id-U_2)  \\
\vdots & & \ddots \
 \ep
\]
\end{theorem}

This formula makes it possible to compute enough $\rho$--invariants to show that 
several interesting boundary links are neither
boundary link slice nor slice. 
Note that if $L$ is boundary link slice then $\rho(M_L,\v)(\a)$ for all representations
$\a$ with $\det(M(A,\a))=0$, i.e. not only for representations that factor through a $p$--group.
Levine announced a proof that this result also holds in the case that $L$ is slice.

The structure of this paper is as follows. In section \ref{sectiontop} we'll give a more detailed exposition
of the eta-invariant and the rho-invariant. In particular we'll cite a criterion of Levine's when homology $G$--bordant
manifolds have identical eta-invariants. These results will be applied in section \ref{sectionetaforlink}
to link concordance questions and in section \ref{sectionboundary} to boundary link concordance questions.
We furthermore define a useful signature function for boundary links. We apply our invariants
to several interesting cases in section \ref{sectionexamples}.
We conclude the paper with two sections containing a formula relating eta-invariants of finite
covers and  the computation of the
$\rho$-invariant for boundary links.

{ \bf Acknowledgment.} I would like to thank Jerry Levine, Desmond Sheiham and Jae Choon Cha for many helpful discussions.

\section{The eta invariant as cobordism invariant} \label{sectiontop}
Let $M^{2q+1}$ be a closed odd-dimensional smooth manifold and $\a:\pi_1(M) \to U(k)$ a unitary representation.
Atiyah, Patodi, Singer \cite{APS75} associated to $(M,\a)$ a number $\eta_{\a}(M)$ called the (reduced)
eta-invariant of
$(M,\a)$. For more details cf. section \ref{sectionetas}.
 
For a hermitian matrix or form $A$ (i.e. $\bar{A}^t=A$) we define
\[ \sign(A):=\# \mbox{ positive eigenvalues of $A$}-\# \mbox{ negative eigenvalues of $A$}\]
and for a skew-hermitian matrix $A$ (i.e. $\bar{A}^t=-A$) we define
$\sign(A):=\sign(iA)$.
 
The main theorem to compute the eta-invariant is the following (cf. \cite{APS75}).
\begin{theorem} (Atiyah-Patodi-Singer index theorem) \label{apsthm}
Let $(M^{2q+1},\a)$ as above.
If there exists $(W^{2q+2},\b:\pi_1(W)\to U(k))$ with $\partial(W^{2q+2},\b)=r(M^{2q+1},\a)$
for some $r\in \N$, then
\[ \eta_{\a}(M) =\frac{1}{r}(\sign_{\b}(W)-k\sign(W)) \]
\end{theorem}

Let $G$ be a group, then a $G$-manifold is a pair $(M,\varphi)$ where $M$ is a compact oriented manifold with
components $\{M_i\}$ and $\v$ is a collection of homomorphisms $\varphi_i:\pi_1(M_i) \to G$ where 
each $\varphi_i$ is
defined up to inner automorphism.
Let $R_k(G):=\{ \a| G \to U(k)\}$.
For an
odd-dimensional $G$-manifold $(M,\varphi)$   define
\[ \ba{rcl} \rho(M,\varphi):{R}_k(G) &\to & \R \\
       \a &\mapsto &\eta_{\a\circ \varphi}(M) \ea \]

  We call two odd-dimensional $G$-manifolds $(M_j,\a_j), j=1,2$,
homology $G$-bordant if there exists a
$G$-manifold
$(N,\b)$  such that $\partial(N)=M_1 \cup -M_2, H_*(N,M_j)=0$ for $j=1,2$
and, up to inner automorphisms of $G$,
$\beta|\pi_1(M_j)=\a_j$.
We want to relate the $\rho$-function for homology $G$-bordant manifolds.

 
Let 
\[ P_k(G)=\{ \a \in R_k(G) | \a \mbox { factors through a group of prime power order}\} \]

\begin{theorem} \cite[p. 95]{L94} \label{levinethm}
If  $(M_i,\a_i), i=1,2$ are homology $G$-bordant manifolds, then 
\[ \rho(M_1,\v_1)(\a)=\rho(M_2,\v_2)(\a) \mbox{ for all }\a \in {P}_k(G)  \]
 \end{theorem}

\section{Eta invariants as link concordance invariants} \label{sectionetaforlink}
Let $L\subset S^{2q+1}$ be a link.
We'll study the eta-invariants associated to the closed manifold $M_L$, the result of zero-framed surgery
along $L\subset S^{2q+1}$.
We first compute the eta invariants of the trivial link.

\begin{lemma} \label{lemmaetazero}
Let $M_O$ be the zero-framed surgery on the trivial link $L$.
Then for any $\a:\pi_1(M_O) \to U(k)$ we
get $\eta_{\a}(M_O)=0$.
\end{lemma}

\begin{proof}
Let $\a:\pi_1(M_O)\to U(k)$ be a representation. Let $D_1,\dots,D_m$ be the push-in off the 
disks in $S^{2q+1}$ bounding $L_1,\dots,L_m$ and let $W:=D^{2q+2}\sm (N(D_1)\cup \dots\cup N(D_m))$.
Note that $\pi_1(S^{2q+1}\sm L)\cong \pi_1(W)\cong F$,
in particular we can use $W$ to compute $\eta_{\a}(M_O)$.
But $W$ is homotopy equivalent to the wedge of $m$ circles,
in particular 
$H_{q+1}(W)=H_{q+1}^{\a}(W, \C^k)=0$,
hence the untwisted and twisted signatures vanish, hence $\eta_{\a}(M_O)=0$ by theorem
\ref{apsthm}.

\end{proof}

\subsection{Abelian eta invariants} \label{sectiongzm}
Recall that any oriented link $L$ with $m$ components has a canonical map $\eps_L: \pi_1(M_L)\to H_1(M_L)=\Z^m$.
Furthermore if $L_1,L_2$ are link concordant, then 
$(M_{L_1},\eps)$ and $(M_{L_1},\eps)$ are canonically homology $\Z^m$-bordant.

The following is now immediate from theorem \ref{levinethm}.

\begin{proposition}  
Let $L_1,L_2$ be concordant links, then  
\[ \rho(M_{L_1},\eps)(\a)=\rho(M_{L_2},\eps)(\a) \mbox{ for all }\a \in P_k(\Z^m) \]
\end{proposition}

The following corollary contains basically  the statement of Smolinsky's main theorem in \cite{S89b}. 
It follows immediately from the proposition and lemma \ref{lemmaetazero}.

\begin{corollary}
Let $L$ be a slice link, $\a \in P_1(\Z^m)$, then $\eta_{\a}(M_L)=0$.
\end{corollary}


\begin{remark}
Levine \cite{L94} shows that that there are links whose eta-invariants vanish for all $\a \in P_1(\Z^m)$
but where a close study of $\rho(M,\eps):R_1(\Z^2)\to \R$ still shows that the
links are not slice.
\end{remark}

We quickly recall a result from high-dimensional knot theory.
Combining results of Matumuto \cite{M77} and Levine \cite{L69}, \cite{L69b}  we get the following theorem.

\begin{theorem} \label{knotu1representations}
If $q>1$, then a knot $K\subset S^{2q+1}$ represents a torsion element in $C(2q-1,1)$
if and only if  $\eta_{\a}(M_K)=0$ for all $\a\in P_1(\Z)$.  
\end{theorem} 

In section \ref{sectionexamples} we show that one--dimensional eta--invariants are not enough 
to detect non--torsion elements in $C_{2q-1}(B_m)$ for $m>1$ and $q>1$.


\subsection{Non--abelian eta invariants} \label{sectiongffi}

For $G$ a group define the lower central series inductively by $G_0:=G$ and $G_i:=[G,G_{i-1}], i>0$.
Milnor \cite{M57} showed that for a link $L$
\[ \pi_1(S^3\sm L)/ \pi_1(S^3\sm L)_k\cong \langle x_1,\dots,x_m | [x_i,w_i], \langle x_1,\dots,x_m\rangle_k \rangle \]
where $x_i$ are representatives for the meridians, $w_i$ for the longitudes and $\langle x_1,\dots,x_m\rangle_k$
denotes the $k^{th}$ term in the lower central series of   the free group generated by
$x_1,\dots,x_m$.

To avoid confusion we'll henceforth denote  the free group on $m$ generators $t_1,\dots,t_m$ by $F$.
Let $F\to \pi_1(S^{2q+1}\sm L)=:\pi$ be a map $t_i$ to a meridian of 
the $i^{th}$ component of $L$. Levine \cite{L94} shows that this induces isomorphisms
$F/F_i \xrightarrow{\cong} \pi/\pi_i$  for all $i$
if $q>1$.
If $q=1$, then
we say that  $L$ has zero $\bar{\mu}$--invariant of level $i$ if this induces
an isomorphism $F/F_i \xrightarrow{\cong} \pi/\pi_i$. By Milnor's result on $\pi_1(S^3\sm L)/ \pi_1(S^3\sm L)_k$
a knot has zero $\bar{\mu}$--invariant of level $i$ if and only if for longitudes
$\l_1,\dots,\l_m$ , $\{ \l_j \} \in \pi_1(S^{2q+1}\sm L)_i$.
Examples for $1$-dimensional links with zero $\bar{\mu}$--invariants are boundary links.

We say $\varphi:\pi_1(S^{n+2}\sm L)\to F/F_i$ is an $F/F_i$--structure   if a 
meridian of the $j^{th}$ component gets 
sent to  $t_j$.  Note that it follows from Stalling's theorem \cite{S65} that conjugates of generators
for $F/F_i$ are also generators of $F/F_i$.

The case $i=1$ is of course uninteresting since $F/F_i=\Z^m$.
If $i>1$  then $L$ has in general no canonical $F/F_i$--structure.

\begin{lemma}\cite[p. 101]{L94} \label{lemmapsi}
If $\varphi_1$ and $\varphi_2$ are
$F/F_i$--structures for the same link, then $\varphi_1=\psi \circ \varphi_2$ for an automorphism of $F/F_i$ that
sends
$t_j$ to a conjugate of
$t_j, j=1,\dots,m$
\end{lemma}  
We call such an automorphism  a special automorphism of $F/F_i$.
A link $L$ equipped with an $F/F_i$--structure is called $F/F_i$--link.
Let $(L_1,\varphi_1),(L_2,\varphi_2)$ be two  $F/F_i$--links, we say they are $F/F_i$--concordant if there exists
a link concordance $C$ and a map $\varphi:\pi_1(S^{2q+1}\times [0,1]\sm C)\to F/F_i$ which restricts to  
$\varphi_1$ and $\varphi_2$ up to inner automorphism.

The following proposition is well--known.
\begin{proposition} \label{proplinkgconc}
\bn
\item
 If  $L_1$ is an $F/F_i$--link and $L_2$ is link concordant to $L_1$,
then there exists an $F/F_i$--structure on $L_2$ such that $L_1$ and $L_2$ are $F/F_i$--concordant.
\item If $L_1,L_2$ are link concordant and $L_1$ has zero $\bar{\mu}$--invariants of level $j$, then 
$L_2$ also has zero $\bar{\mu}$--invariants of level $j$.
\item A one-dimensional slice link has zero $\bar{\mu}$--invariant for all levels.
\en
\end{proposition}

\begin{proof}
Let  $C\subset S^{2q+1}\times [0,1]$ be a link concordance between $L_1$ and $L_2$.
\bn
 \item Consider 
\[ \pi^j:=\pi_1(S^{2q+1}\sm L_j)\to \pi_1(S^{2q+1}\times [0,1]\sm C)=:\pi_C \]
These maps are normally surjective  and hence define isomorphisms
$\pi_C/\pi_{C,i}\cong \pi^j/\pi^j_i\cong F/F_i$ by Stalling's theorem \cite{S65}. The statement now follows easily
(cf. \cite[p. 102]{L94} for details).  
\item 
This follows immediately from the definition and $ F/F_i\cong\pi^1/\pi^1_i\cong \pi_C/\pi_{C,i}\cong \pi^2/\pi^2_i$.
\item This follows immediately from $(2)$ since a slice link is concordant to the unlink which has obviously
zero $\bar{\mu}$-invariant for all levels.
\en
\end{proof}

It is clear that in the case $q>1$ the map $\pi_1(S^{2q+1} \sm L)\to \pi_1(M_L)$ is an isomorphism,
hence
\[ \pi_1(M_L)/\pi_1(M_L)_i \cong \pi_1(S^{2q+1}\sm L)/\pi_1(S^{2q+1}\sm L)_i  \]
If $q=1$
the kernel $\pi_1(S^3 \sm L)\to \pi_1(M_L)$ is generated by the longitudes.
In particular if $L$ has zero $\bar{\mu}$-invariants of level $i$, then
\[ \pi_1(M_L)/\pi_1(M_L)_i \cong \pi_1(S^3\sm L)/\pi_1(S^3\sm L)_i  \]
 In both cases an $F/F_i$--structure on $L$ gives an $F/F_i$--structure on $M_L$.

\begin{proposition}\cite[p. 102]{L94} \label{propmlibordant}
If $(L_1,\v_1), (L_2,\v_2)$ are $F/F_i$--concordant $F/F_i$--links, then $(M_{L_1},\v_1)$ and $(M_{L_2},\v_2)$ are homology
$F/F_i$--bordant.
\end{proposition}

\begin{proof}
If $C$ is an $F/F_i$--concordance, then doing surgery along $C\subset S^{2q+1}\times [0,1]$ gives 
a homology $F/F_i$--bordism for $(M_{L_1},\v_1)$ and $(M_{L_2},\v_2)$.
\end{proof}

 
The following is immediate from theorem \ref{levinethm}, lemma \ref{lemmapsi}  and propositions \ref{proplinkgconc},
\ref{propmlibordant}. The theorem generalizes results on link concordance by Cha and Ko \cite{CK99}.

\begin{theorem} \label{thmffilinks}
Let $L_1,L_2$ be concordant links. If $\varphi_1,\varphi_2$ are arbitrary $F/F_i$--structures for $L_1,L_2$, 
then there exists a special automorphism $\psi$ of $F/F_i$ such that
\[ \rho(M_{L_1},\v_1)(\a)=\rho(M_{L_2},\psi \circ \v_2 )(\a) \mbox{ for all }\a \in P_k(F/F_i) \]
\end{theorem}

\subsection{Representations of $F/F_2$} 
We now give an example of a non-trivial (i.e. non-abelian) unitary   representation of $F/F_2$.
For $U_1,\dots,U_m  \in U(k)$ define $ \a_{(U_1,\dots,U_m)}: F\to U(k)$ by $\a(t_i):=U_i$. We'll find $U_1,\dots,U_m$
such that $\a_{(U_1,\dots,U_m)}$ factors through $F/F_2$.

Let $z_1,\dots,z_k\in S^1$ and $\chi:F\to S^1$ a character such that $\chi(t_i^k)=1$.
Define 
\[ 
 U_1:= \bp 0&\dots &0&z_k \\
z_1&\dots &0&0\\
 0&\ddots&&\vdots \\
0&\dots&z_{k-1}&0\ep, \quad
U_i:=\bp \chi(t_i)&0&\dots&0 \\
0&\chi(t_1t_i)&&0\\
0&0&\ddots &\vdots \\
0&0&\dots &\chi(t_1^{k-1}t_i)\ep, i=2,\dots,m
\]
\begin{lemma} \label{lemmaff2rep}
The representation $\a=\a_{(U_1,\dots,U_m)}: F\to U(k)$   factors through $F/F_2$.
\end{lemma}

\begin{proof}
It is clear that we are done once we show that for all $x\in [F,F]$, $\a(x)\in \C \cdot \id$.
Since
 \[ [x,vw]=[x,v]v[x,w]v^{-1}\]
we only have to show that $\a([x_i,x_j])\in \C\cdot \id$,
but an easy calculation using $\chi(t_i^k)=1$ shows that
\[ \ba{rcll} \a([t_1,t_j])&=&\chi(t_1^{-1})\cdot \id & \mbox{if } j\ne 1 \\
\a([t_j,t_1])&=&\chi(t_1)\cdot \id & \mbox{if } j\ne 1 \\
\a([t_i,t_j])&=&\id & \mbox{if } i\ne 1 \mbox{ and }j\ne 1\ea \]
\end{proof}

Let $p$ a prime, $k$ a power of $p$, and $z_1,\dots,z_k,\chi$ such that $z_1^{p^N}=\dots=z_k^{p^N}=1$  and
$\chi(v)^{p^N}=\id $ for some $N$, then $\varphi\in P_k(F/F_2)$.
Such a representation turns out to discover non-slice knots in many interesting cases.

This example can easily be generalized to give more complex representations of $F/F_2$.



\subsection{Sliceness obstructions} 

\begin{theorem} \label{theoremfirstslice}
Let $L\subset S^{2q+1}$ be a slice link and let $\a \in P_k(\pi_1(M_L))$, then $\eta_{\a}(M_L)=0$.
  \end{theorem}

\begin{proof}
Assume that $\a$ factors through a $p$-group $P$. Then $P_i=\{e\}$ for some $i$  since any $p$-group
is nilpotent  (cf. \cite[p. 169]{J97}). In particular $\a$ factors through $\pi_1(M_L)/\pi_1(M_L)_i$
which is isomorphic to $F/F_i$ since any slice link has zero $\bar{\mu}$--invariants
by proposition \ref{proplinkgconc}. Henc $\a=\b\circ \v$ for some $F/F_i$--structure $\v$ and
some representation $\b$. 
The statement now follows
immediately from proposition \ref{proplinkgconc}, theorem  \ref{thmffilinks} and lemma \ref{lemmaetazero}
since a slice link is concordant to the unlink.
\end{proof}

Define $PD(k)\subset U(k)$ to be the subgroup generated by permutation matrices and diagonal matrices.
For a prime $p$ define $PD_p(k)\subset PD(k)$ to be the subgroup of matrices where all eigenvalues are roots of unity of order
a power of $p$. It is generated by all permutation matrices whose order is a power of $p$ and all diagonal matrices
whose entries are roots of unity of order
a power of $p$. Note that a finitely generated subgroup $PD_p(k)$ is in fact a finite group, hence a $p$-group.


\begin{theorem} \label{theoremffius}
Let $L\subset S^{2q+1}$ be a slice link with meridians $\mu_1,\dots,\mu_m$. Let $p$ be a prime number
and let $U_1,\dots,U_m\in PD_p(K)$. Then there exists a unique representation $\b:\pi_1(M_L)\to U(k)$ with
$\b(\mu_j)=U_j$. Furthermore $\eta_{\b}(M_L)=0$.
\end{theorem}

\begin{proof}
Let $\a:=\a(U_1,\dots,U_m):F\to U(k)$, then $\im(\a)$ is a $p$-group, hence $\a$ factors through $F/F_i$ for some $i$.
It's clear that $\b$ is given by $\pi_1(M_L)/\pi_1(M_L)_i\cong F/F_i\to U(k)$. Furthermore $\b\in P_k(\pi_1(M_L))$,
the theorem now follows from theorem \ref{theoremfirstslice}.
\end{proof}

\begin{proposition}
  Let $\a \in P_k(F/F_i)$, then there exists a prime $p$ 
such that $\a$ is conjugate to a representation $\ti{\a}$ with $\ti{\a}(t_j)\in
PD_p(k)$ for all $j$.
 \end{proposition}

\begin{proof}
  This follows from the fact that if $\a:P\to U(k)$ is a representation of a $p$-group $P$, then $\a$ is  
induced from a representation of degree 1 (cf. \cite[p. 578ff]{H67}). This means that there exists a subgroup $Q\subset P$
and a one-dimensional representation $Q\to U(\C)$ such that $\a$ is given by
the natural $P$-left action on $\C P\otimes_{\C Q} \C$. Pick representatives $p_1,\dots,p_k$ for
$P/Q$, writing $\a$ with respect to  this basis we see that $\a$ is of the required type.
 \end{proof}


\begin{remark}
The above proposition together with theorem \ref{theoremfirstslice} shows that theorem \ref{theoremffius} is the best
possible sliceness obstruction theorem which can be based on theorem \ref{levinethm}. 
\end{remark}

\subsection{Algebraic closures of groups and link concordance}
Whereas theorem \ref{theoremffius} can't be improved on with our means there's still room for improvement for 
proposition \ref{thmffilinks} because of the extra indeterminacy given by the special automorphism group.

For a group $G$ Levine \cite{L89a}, \cite{L89b}, \cite{L90} introduced the notion of algebraic closure $\hat{G}$
and residually nilpotent algebraic closure $\bar{G}$ of a group $G$.
The results of section \ref{sectionetaforlink} for $G=F/F_i$ also hold for $G=\hat{F}$ and $G=\bar{F}$
(cf. \cite{L94} for details), in particular links with zero $\bar{\mu}$-invariants  
have a $\bar{F}$-structure and   concordant links are also $\bar{F}$-concordant, same for $\hat{F}$.
In particular we get link concordance invariants from representations in $P_k(\bar{F})$ and $P_k(\hat{F})$.

Note that $p$-groups are nilpotent and hence its own algebraic closure
(\cite[p. 100]{L90}). This shows that representations in $P_k(\pi_1(M_K))$ that factor through an $F/F_i$--structure for some
$i$ correspond to representations that factor through some $\bar{F}$--structure (or $\hat{F}$--structure).

The following theorem is a stronger version of \ref{thmffilinks}
 
\begin{theorem}
Let $L_1,L_2$ be concordant links with vanishing $\bar{\mu}$--invariants. If $\varphi_1,\varphi_2$ are arbitrary
$\bar{F}$--structures for
$L_1,L_2$,  then there exists a special automorphism $\psi$ of $\bar{F}$ such that
\[ \rho(M_{L_1},\v_1)(\a)=\rho(M_{L_2},\psi \circ \v_2 )(\a) \mbox{ for all }\a \in P_k(\bar{F}) \]
\end{theorem}

\subsection{Relation to previous link concordance invariants}
One can easily see that theorem \ref{theoremffius} contains the sliceness obstructions defined
by Smolinsky   \cite{S89}, \cite{S89b}.

We quickly recall a results by Cha and Ko   and show how it follows from our results.
 
\begin{theorem}\cite[thm. 7]{CK99} \label{thmchako}
Let $L$ be a slice link and $p$ a prime. Let $\varphi:\pi_1(M_L)\to G$ be a homorphism to a finite abelian $p$-group $G$.
Denote the $G$-fold cover of $M_L$ by $M_G$. Let $\a_G:H_1(M_G)\to \Z/p\to U(1)$ be a representation that factors
through
$\Z$, then 
\[ \eta(M_G,\a_G)=0 \] 
\end{theorem}

\begin{proposition}
If a link $L$ satisfies the conclusion of theorem \ref{theoremfirstslice} then it also satisfies the conclusion of theorem
\ref{thmchako}.
\end{proposition}

\begin{proof}
Let $s=|G|$. 
By theorem \ref{thmetamprimem}
there exists a   unitary representation $\a:\pi_1(M_L)\to U(s)$ such that
 \[ \eta_{\a_G}(M_G)=\eta_{\a}(M_L)-s \eta_{\a(G) }(M_L) \]
where $\a(G)$ stands for the representation $\pi_1(M_L)\to U(\C[\pi_1(M_L)/\pi_1(M_G)])=U(\C G)$ given by left
multiplication.
Furthermore $\a \in P_{s}(\pi_1(M))$ by lemma \ref{lemmaalphaprimetoo} and $\a(G)\in P_1(\pi_1(M))$ since $G$ is of prime power
order.

If a link $L$ satisfies the conclusion of theorem \ref{theoremfirstslice},  then $\eta_{\a(G)}(M_L)=0$ and $\eta_{\a}(M_L)=0$.
\end{proof}

In later, unpublished work Cha showed that if $L$ is a slice link, $p$ a prime power, $M'$ a $p^a$--cover of $M_L$ 
(not necessarily regular) and $\a':H_1(M')\to U(1)$ a character whose order is a power of $p$, then $\eta(M',\a')=0$.
In this case we can find $M_L=M_0\subset M_1\subset \dots\subset M_k=M'$ such that $M_i/M_{i-1}$ is a regular
$p$--covering. Using lemma \ref{lemmaalphaprimetoo} and theorem  \ref{thmetamprimem} one can inductively write
$\eta(M',\a)$ as a sum of eta invariants of $M_L$ with representations factoring through $p$--groups.
This shows that Cha's extended result is contained in theorem \ref{theoremfirstslice}.

\section{Eta-invariants and signatures of boundary links} \label{sectionboundary}
 
\subsection{Eta-invariants as boundary link concordance invariants} \label{sectionetaforboundaryslice}
In this section we denote the free group on $m$ generators once again by $F_m$.
Let $(L,\v)\subset S^{2q+1}$ be an $F_m$--link. If $q>1$ then $\pi_1(S^{2q+1}\sm L)\to \pi_1(M_L)$ is an isomorphism.
If $q=1$, then $\v(\l)=e$ for any longitude, since
$[\l_i,\mu_i]=1 \in \pi_1(S^3\sm L)$.
In particular for any $q$ the map $\v$ factors through
$\pi_1(M_L)$.

\begin{proposition}\cite[p. 102]{L94}
Let $(L_1,\v_1), (L_2,\v_2)$ be $F_m$-concordant links, then $(M_{L_1},\v_1), (M_{L_2},\v_2)$ are homology 
$F_m$-bordant.
\end{proposition}

The following theorem is immediate from lemma \ref{lemmacam}, proposition \ref{propa2trivial}, theorem \ref{levinethm}
and the above proposition.

\begin{theorem} \label{thmbdyslice}
Let $(L_1,\varphi_1)$ and  $(L_2,\varphi_2)$ be $F_m$-concordant $F_m$-links, then 
\[ \rho(M_{L_1},\v_1)(\a)=\rho(M_{L_2},\v_2)(\a) \mbox{ for all }\a \in P_k(F_m) \]
If $L_1,L_2$ are boundary concordant boundary links with {\it two} components, then 
\[ \rho(M_{L_1},\v_1)(\a)=\rho(M_{L_2},\v_2)(\a) \mbox{ for all }\a \in P_k(F_2) \]
for any $F_2$-structures $\v_1$ and $\v_2$.
\end{theorem}

The following is immediate from theorem \ref{theoremfirstslice}.

\begin{theorem}\label{thmetaforbdylinkslice}
If $L$ is a boundary link, and $L$ is slice (in particular if $L$ is boundary slice),
then
\[ \rho(M_L,\v)(\a)=0 \mbox{ for any }\a \in P_k(F_m)\]
for any $F_m$-structure $\v$.
\end{theorem}

\begin{corollary}
If $L_1,L_2$ are boundary link concordant boundary links and if $\v_1,\v_2$ are $F_m$-structures, then
there exists a special automorphism $\psi\in CA_m$ such that
\[ \rho(M_{L_1},\v_1)(\a)= \rho(M_{L_1},\psi \circ \v_1)(\a) \mbox{ for any }\a \in P_k(F_m) \]
\end{corollary}

\begin{proof}
Levine \cite[p. 102]{L94} showed that 
if $L_1$ is an $F_m$-link and  $L_2$ a boundary link which is boundary link concordant
to $L_1$, then there exists an $F_m$-structure on $L_2$ such that $L_1$ and $L_2$ are $F_m$-concordant.
The corollary now follows from lemma  \ref{lemmacam} and theorem \ref{thmbdyslice}.
\end{proof}


In section \ref{sectionetaforboundary} we compute the $\rho$-invariant for an $F_m$-link.
This will involve the computation of the eta-invariant of a circle which necessitates the definition
of the following function.
Let $z=e^{2\pi ia}\in S^1$ with $a\in [0,1)$, then define
\[ \eta(z):=\left\{ \ba{rl} 0 &\mbox{ if }a=0 \\ 1-2a &\mbox{ if }a>0 \ea \right. \]
Now we can formulate the following theorem which will be proven in section \ref{sectionetaforboundary}.

\begin{theorem} \label{thmbdyetacomp}
Let $(L\subset S^{2q+1},\v)$ be an $F_m$-link, 
$A=(A_{ij})_{i,j=1,\dots,m}$ a Seifert matrix for $(L,\v)$ of size $(g_1,\dots,g_m)$,
 $\a:F_m\to U(k)$ a representation.
Let  $\eps:=(-1)^{q+1}, g:=\sum_{i=1}^m g_i, T:=\diag(t_1,\dots,t_1,\dots,t_m,\dots,t_m)$ where each $t_i$ appears $2g_i$
times. Let
$\{z_{ij}\}_{j=1,\dots,k}$ be the set of eigenvalues of
$\a(t_i)$. Then 
\[ \ba{rcccl} \rho(M_L,\v)(\a)&=&
&\eps\sum_{i=1}^m sign(\sqrt{\eps}(A_{ii}+\eps A_{ii}^t))\sum_{i=1}^m \sum_{j=1}^k \eta(z_{ij})
+\\
&&+&\sign( \sqrt{-\eps}( A -\eps  \a(T) A^t\a(T)^{-1} -A\a(T)^{-1}   + \eps  \a(T)A^{t})) \ea \]
 where we consider $A$ as a $2gk \times 2gk$ matrix, where each entry of $A=(a_{ij})$ is 
replaced by $a_{ij} \cdot \id_{k}$.
This simplifies for $\eps=-1$ to the following
\[ \rho(M_L,\v)(\a)=
\sign( A + \a(T) A^t\a(T)^{-1} -A\a(T)^{-1}-  \a(T)A^{t}) \]
\end{theorem}

Note that if we let $U_i:=\a(t_i)$, then $  A -\eps  \a(T) A^t\a(T)^{-1} -A\a(T)^{-1}   + \eps  \a(T)A^{t}$ equals 
\[
\bp A_{11}(1-U_1^{-1})-\eps A_{11}^t(1-U_1 ) & A_{12}(1-U_1)(1-U_2^{-1}) & \dots \\
A_{21}(1-U_2)(1-U_1^{-1}) &A_{22}(1-U_2^{-1})-\eps A_{22}^t(1-U_2)  \\
\vdots & & \ddots \
 \ep
\]
here we use the convention of the theorem again, i.e. we view $A_{ij}$ as a $2g_ik\times 2g_jk$--matrix.
Alternatively we could write $A_{11}\otimes(1-U_1^{-1})-\eps A_{11}^t\otimes (1-U_1)$ etc..
 
This result generalizes a computation done by Cha and Ko \cite{CK99} for certain unitary representations.
We suggest the following conjecture which would be a generalization of theorem \ref{knotu1representations}.

\begin{conjecture}
Let $q>1$ and $(L,\v)\subset C_{2q+1}(F_m)$. 
If for all $k$, $\rho(M_L,\v)(\a)=0$ for a dense set of
representations $\a \in R_k(F_m)$, then $(L,\v)$ represents a torsion element.
\end{conjecture}

Note that in light of theorem  \ref{theoremboundaryclass} this conjecture is purely algebraic.
This conjecture seems to be hard to prove, and any attempt would I think require a good understanding
of non-commutative algebraic geometry. 
If this conjecture can be proven to be true then this would give an   algorithm for detecting
non--torsion elements in $C_n(F_m)$, which is easier to implement
than Sheiham's \cite{S02} algorithm. The disadvantage of 
such an algorithm would be that it can not conclude   in finite time that an $F_m$--link is torsion.  
An interesting follow up
question to a positive answer would be whether there exists a $k$ depending computably on $(L,\v)$, 
such that it is enough to
study the
$\rho$--invariant for dimensions less or equal than $k$ for deciding whether $(L,\v)$ is torsion or not.

\subsection{Signature invariants for boundary link matrices}
Recall that if a boundary link $(L,V)$ is boundary link slice  
then any Seifert matrix is metabolic. Using this fact and some algebra we can strengthen
theorem \ref{thmetaforbdylinkslice}.

Let $A=(A_{ij})$ be an $\eps$-Seifert matrix and $U_i\in U(k), i=1,\dots,m$.
We denote by $U:=\diag(U_1,\dots,U_m)$ the block diagonal
matrix with blocks
$U_i\cdot \id_{h_i}$ and
define 
\[ M(A,U):=\sqrt{-\eps}( A -\eps UA^tU^{-1}-AU^{-1}+\eps  UA^{t})\]
using the convention of theorem \ref{thmbdyetacomp}. Furthermore let $\s(A,U):=\sign(M(A,U))$. 
If $A$ is metabolic then $M(A,U)$ is metabolic as well, if $U$ is such that $\det(M(A,U))\ne 0$ then $\s(A,U)=0$.
The map $\s$ is continuous outside of the set
\[ S_k(A):=\{ (U_1,\dots,U_m) \in U(k)^m | \det(M(A,U))=0\}   \]
It is easy to see that if $A_1,A_2$ are S-equivalent, then $\s(A_1,U)=\s(A_2,U)$ and 
$S(A_1)=S(A_2)$. In particular
for a boundary link pair $(L,V)$ we can define $\s(L,V,U):=\s(A,U)$ 
using any Seifert matrix and we let $S_k(L,V):=S_k(A)$.
This generalizes signature invariants for knots defined by 
Levine \cite{L69} and Trotter \cite{T73}. 

We immediately get the following proposition.
\begin{proposition} \label{propsignbdyslicezero}
Let $(L,V)$ be a boundary link pair which represents zero in $C_n(B_m)$, then 
$\s(L,V,(U_1,\dots,U_m))=0$ for all $(U_1,\dots,U_m) \not\in S_k(L,V)$. 
\end{proposition}

Combining this with theorem 
\ref{thmbdyetacomp} we get a theorem that gives   a much stronger boundary sliceness obstruction
than theorem \ref{thmetaforbdylinkslice} since the matrices $U_i$ no longer have to lie in $PD_p(k)$ for some prime $p$.

\begin{theorem} \label{thmbdylinkzero}
Let $(L,V)$ be a boundary link pair which represents zero in $C_n(B_m)$, then 
\[ \rho(M_L,\v)(\a_{(U_1,\dots,U_m)})=0 \mbox{ for all }(U_1,\dots,U_m) \not\in S(L,V) \] 
\end{theorem}

\begin{remark}
Note that $\sqrt{-\eps}( A -\eps UA^tU^{-1}-AU^{-1}+\eps  UA^{t})=\sqrt{-\eps}(A+\eps UA^t)(1-U^{-1})$, using an argument
as in \cite{L69} one can show in a purely algebraic way that $\s$ is continuous outside of the set 
$\{ (U_1,\dots,U_m) \in U(k)^m | \det(A+\eps UA^t )=0\} $. This agrees with the topological result of Levine's
\cite{L94} that $\rho$ is continuous on the
set
\[ \{\a=\a_{(U_1,\dots,U_m)}| H_1^{\a}(M_L,\C^k)=0 \}\]
since $A+\eps TA^t$ represents the homoloy of the universal $F_m$--cover 
of $M_L$. This shows in particular that $\rho$ is zero in a neighborhood of the trivial representation.
\end{remark} 

\begin{remark}
Levine announced a proof that these $\rho$--invariants are in fact obstructions to a link being slice and not just
obstructions to a link being boundary slice.
\end{remark}

There are many ways to associate a hermitian matrix to $A$ which is metabolic if $A$ is metabolic.
Let 
$F_{ij}\in M(k_i\times k_j,\C), i,j=1,\dots,m$ such that 
$F_{ji}=\sqrt{-\eps} \bar{F}_{ij}^t$, then we also get a similar proposition for
\[ \s(A,F_{ij}):=\sign \bp A_{11}F_{11}+A_{11}^t\bar{F}_{11}^t & A_{12}F_{12}&\dots \\ 
A_{21}F_{21} & A_{22}F_{22}+A_{22}^t\bar{F}_{22}^t \\
\vdots &&\ddots \ep \]
This approach has the advantage that it is much easier to find (random) matrices in $M(k_i\times k_j,\C)$ than
matrices in $U(k)$.

In the knot case one can easily show that these signature invariants have the same information content
as $\s(A,U)$. In the case $m>1$  we don't know
whether these different signature functions have different information content or not.

\section{Examples} \label{sectionexamples}
Ko \cite{K87} gives an example of a three component boundary link $L\subset S^{4l+3}$ 
with Seifert manifold $V$ such that
$(L_{1,-2},V_{1,-2}):=(L,V)\#- (L,\a_{12}V)$ (cf. \cite{K87} for details on the action of $CA_3$ on Seifert surfaces)
  has the following Seifert matrix of size $(1,2,2)$
\[ A=\bp  
0&1 &1&0&0&0& 0&0&0&0\\
0&0 &0&0&0&0& 1&0&0&0\\
1&0 &0&1&0&0& 1&0&0&0\\
0&0 &0&-1&0&0& 0&0&0&0\\
0&0 &0&0&0&-1& 0&0&-1&0\\
0&0 &0&0&0&1& 0&0&0&0\\
0&1 &1&0&0&0& 0&1&0&0\\
0&0 &0&0&0&0& 0&-1&0&0\\
0&0 &0&0&-1&0& 0&0&0&-1\\
0&0 &0&0&0&0& 0&0&0&1\ep \]
Ko showed that $L_{1,-2}$ is not boundary slice and posed the question whether $L_{1,-2}$ is slice or not.
By construction we get for $\a \in P_1(F_3)$ that
\[
\rho(M_{L_{1,-2}},\v_{V_{1,-2}})(\a)=
\rho(M_{L},\v_{V})(\a)-\rho(M_{L},\a_{12}\v_{V})(\a)
\rho(M_{L},\v_{V})(\a)-\rho(M_{L},\v_{V})(\a \circ \a_{12})=0 \]
since $U(1)$ is abelian. Hence all one--dimensional eta--invariants vanish.

Cha and Ko \cite{CK99} showed that $L$ is in fact not slice. We reprove this using higher dimensional representations.
Let
\[ U_1=\bp 0 &1 \\ 1 &0\ep, \quad U_2=\bp \frac{\sqrt{2}}{2}+\frac{\sqrt{2}}{2}i&0\\
0&-\frac{\sqrt{2}}{2}-\frac{\sqrt{2}}{2}i \ep, \quad U_3=
U_3=\bp \frac{\sqrt{2}}{2}+\frac{\sqrt{2}}{2}i&0\\
0&-\frac{\sqrt{2}}{2}-\frac{\sqrt{2}}{2}i \ep
\]
A computation using theorem  \ref{thmbdyetacomp} shows that $\rho(M_L,\v)(\a_{U_1,U_2})=-2$, hence
$L$ is not slice by theorem \ref{theoremffius}.
 
On the other hand,  let $(L_{1,-1},V_{1,-1}):=(L,V)\#- (L,V)$,
then $L_{1,-1}$ is obviously slice. This proves the following proposition.

\begin{proposition} \label{propnogroup}
Connected sum is not a well-defined operation on $C(n,m)$ for $m\geq 3$.
\end{proposition}


We now give an example of a two component link with vanishing one--dimensional eta--invariant but which
is not slice.
Consider the following boundary link Seifert
matrix of size
$(2,1)$:
\[ A=(A_{ij})_{i,j=1,2}=
\bp 0 &0&0&0&0&0\\
1&0&0&-1&0&0\\
0&0&1&0&0&1\\
0&1&-1&-2&0&0\\
1 &0&0&0&-2&0\\
0&0&1&0&1&1\ep \]
 
Let  $(L,V)=(L_1 \cup L_2,V_1 \cup V_2) \subset S^{4l+3}$ be a boundary link  pair
with Seifert matrix $A$. In fact we can find isotopic slice
knots $L_1,L_2$ and corresponding Seifert surfaces with the above property
since one can easily see that $A_{11}$ and $A_{22}$ are $S$-equivalent and metabolic.
 
Note that
$\Delta(L)(t_1,t_2)=\det(AT-A^t)t_1^{-2}t_2^{-1}=-(t_1t_2+t_1^{-1}t_2^{-1})-(t_1^{-1}t_2+t_1t_2^{-1})+5$.
Let $(\ti{L},\ti{V})=(L_2,V_2)\cup (L_1,V_1)$,
clearly $(\ti{L},\ti{V})$ is a boundary link with Seifert matrix
$(\ti{A}_{ij})=(A_{ji})$.

Now pick an arc connecting $L$ and $\ti{L}$ which lies outside of $V$ and $\ti{V}$. Use this arc to
form $L\# -\ti{L}$. If $q>1$ then this link is independent of the choice of the arc.

\begin{proposition} \label{propexampleu1zero}
The boundary link $(L\# -\ti{L},V\# -\ti{V})$ has zero $U(1)$-eta invariants but is not boundary link slice.
Furthermore $L\# -\ti{L}$ is not slice.
 \end{proposition}

\begin{proof}
Let $B=A \oplus -\ti{A}$ be a Seifert matrix for $(L\# -\ti{L},V\# -\ti{V})$.
For $z_1,z_2\in S^1$ let $Z=\diag(z_1,z_1,z_1,z_1,z_2,z_2)$, then
\[ \rho(M_{L\# -\ti{L}},\eps)_{\a_{(z_1,z_2)}}=\sign(B(1-Z)+B^t(1-Z^{-1}))=\sign((BZ-B^t)(Z^{-1}-1)) \]
In particular the function $\rho(M_{L\# -\ti{L}},\eps):R_1(\Z^2)=S^1\times S^1\to \Z$ is constant outside of
the set $\{ (z_1,z_2)\in S^1\times S^1 | \det(AZ-A^t)=0\}$.
It is obvious that for all $z_1,z_2\in S^1$ we have 
\[ \det(BZ-B^t)z_1^{-2}z_2^{-1}=(-(z_1z_2+z_1^{-1}z_2^{-1})-(z_1^{-1}z_2+z_1z_2^{-1})+5)^2 \geq 1 \]
hence the $\rho$-invariant function is constant. Picking $z_1=-1,z_2=-1$ we can compute that
the constant is in fact $0$.

Now let 
\[ U_1=\bp 0 &i \\ 1 &0\ep, \quad U_2=\bp i&0\\
0&-i \ep
\]
A computation using theorem  \ref{thmbdyetacomp} shows that $\rho(M_{L\# -\ti{L}},\v)(\a_{U_1,U_2})=-2$, hence
$L\# -\ti{L}$ is not slice by theorem \ref{theoremffius}.
\end{proof}


Now consider the following Seifert matrix of size $(1,1)$
\[ A= \bp 0&1&0&1 \\ 0&0&2&1 \\ 0&2&1&0 \\ 1&1&1&0 \ep\]  
Let $(L,V)$ be a boundary link pair with Seifert matrix $A$.
If we let 
\[ F_{11}=\bp 4&1\\ 3&0 \ep, \quad F_{12}=\bp 4&1\\ 2&1 \ep, \quad F_{22}=\bp 1&2\\ 4&1 \ep \]
then $\s(A,F_{ij})=-2$, which shows that $A$ is not metabolic, i.e. $L$ is not boundary slice.
 Computer computations indicate that $\rho(M_L,\v)$ vanishes on 
$S_1(L,V)$ and $S_2(L,V)$ but is non-zero on $S_3(L,V)$ which shows again that $L$ is not boundary
slice by proposition \ref{propsignbdyslicezero}.
 
All the $\rho$-invariants of theorem \ref{theoremffius}, i.e. all eta invariants corresponding to representations
that factor through a $p$--group  that I
computed so far with a computer vanish. So it seems like one can not use theorem \ref{theoremffius} to
say that $L$ is not slice. 

A new result by Levine (cf. the second remark after theorem  \ref{thmbdylinkzero}) shows that $L$ is in fact
not slice.


\section{Relating eta-invariants of finite covers}  \label{sectionfinitecovers} \label{sectionetas}
Let $M$ be an oriented Riemannian manifold of dimension $2l-1$ and $\a:\pi_1(M)\to U(k)$ a representation.
Denote the universal cover of $M$ by $\ti{M}$.
Then let $V_{\a}:=\ti{M}\times_{\pi_1(M)} \C^k$, this is a $\C^k$-bundle over $M$.
On the space of differential forms of even degree there's a natural self-adjoint operator $B$ defined
by
\[ \ba{rcl} \Omega_{2k}(M)&\to &\Omega_{2l-2k}(M) \\ 
\w &\mapsto& i^l(-1)^{k+1}(*d-d*)\w \ea \]
This can be naturally extended to give a self-adjoint operator $B_{\a}$ acting on even forms 
with coefficients in the flat vector bundle defined by $\a$. 
Consider the spectral function $\eta_{\a}(M,s)$ of this operator defined by
\[ \eta_{\a}(M,s):=\sum_{\l \ne 0}(\sign(\l))|\l|^{-s} \]
where $\l$ runs over the eigenvalues of $B_{\a}$.
Atiyah-Patodi-Singer \cite{APS75} 
showed that for $s$ with $\mbox{Re}(s)$ big enough, $\eta_{\a}(M,s)$ converges  to a holomorphic function. Furthermore this
holomorphic function can be extended to $0$ and $\eta_{\a}(M,0)$ is finite.
Now define the (reduced) eta--invariant of $(M,\a)$ to be
\[ \eta_{\a}(M):=\eta_{\a}(M,0)-k\eta(M,0)\]
where $\eta(M,s)$ denotes the eta function corresponding to the trivial   one--dimensional
representation of $\pi_1(M)$. Atiyah-Patodi-Singer \cite{APS75} 
showed that   $\eta_{\a}(M)$ is independent of the Riemannian
metric on $M$.\\

Let $M$ be a  manifold of dimension $2l-1$ and $M'$ a finite cover, not necessarily regular.
Let $\a':\pi_1(M')\to U(k)$ be a representation.
The goal is to express $\eta_{\a'}(M')$ in terms of eta-invariants of $M$.

Consider $\C \pi_1(M)\otimes_{\C \pi_1(M')} \C^k$ where we view $\C^k$ as
a $\C \pi_1(M')$-module via $\a'$. 
We give $\C \pi_1(M)\otimes_{\C \pi_1(M')} \C^k$ the metric induced by
\[  ((p_1\otimes v_1),(p_2\otimes v_2))\to \sum_{g\in \pi_1(M')} \delta_{(p_1g,p_2)}\overline{(\a'(g)^{-1}v_1)}^tv_2 
\]
wgere $p_i\in \pi_1(M), v_i\in \C^k$.
It's easy to see that this is well-defined.
Let $s:=[\pi_1(M):\pi_1(M')]$, then clearly  
$\dim(\C \pi_1(M)\otimes_{\C \pi_1(M')} \C^k)=ks$.

Define 
\[ \ba{rcl}
\a:\pi_1(M) &\to& \Aut(\C \pi_1(M)\otimes_{\C \pi_1(M')} \C^k)\\
a &\mapsto &(p\otimes v \mapsto ap \otimes v) \ea
 \]
This action is obviously isometric, i.e. unitary.

Denote by $\a(M,M')$ the representation $\pi_1(M)\to U(\C \pi_1(M)\otimes_{\pi_1(M')}\C)$ given by left
multiplication where we consider $\C$ as the trivial $\pi_1(M')$--module.
\begin{theorem} \label{thmetamprimem}
  \[ \eta_{\a'}(M')=\eta_{\a}(M)-k \eta_{\a(M,M') }(M) \]
\end{theorem}

\begin{proof}
Give $M$ some Riemannian structure and $M'$ the induced structure.
 We have to show that
\[ \eta_{\a'}(M',0)-k\eta(M',0)=(\eta_{\a}(M,0)- ks\eta(M,0))-k (\eta_{\a(M,M')}(M,0)-s\eta(M,0)) \] 
We'll in fact show that
\[ \ba{rcl} \eta_{\a'}(M',0)&=&\eta_{\a}(M,0) \\
\eta(M',0)&=&\eta_{\a(M,M')}(M,0)
\ea \]

Recall that
\[ \ba{rcl}
V_{\a'}&=&\ti{M}\times_{\pi_1(M')} \C^k \mbox{ and }\\
V_{\a}&=&\ti{M}\times_{\pi_1(M)} (\C \pi_1(M)\otimes_{\C\pi_1(M')} \C^k) \ea \]
Let $p\in M, U\subset M$ a (small) neighborhood and $p_1,\dots,p_s,U_1,\dots,U_s$ the  different lifts.
Then the map
\[ \ba{rcl} \oplus_{i=1}^s \pi_i:V_{\a}|_U&\to &\oplus_{i=1}^s V_{\a'}|_{U_i}\\
\sum (q,g_ih_i\otimes v_i)&\mapsto & \sum (qg_i,h_iv_i)\ea \]
is an isomorphism with inverse map given by
\[
\sum (q_i,1\otimes v_i) \mapstoleft \sum (q_i,v_i)
\]
where $g_i$ such $\pi(qg_i)\in U_i$ and $h_i \in \pi_1(M')$.
Note that
\[  
\Omega^{2i}(M',V_{\a})|_{\cup U_i}=\oplus_{i=1}^s \Omega^{2i}(M')|_{U_i}\otimes_{C^{\infty}(U_i)} \Gamma(V_{\a'}|_{U_i})
\]
This is isomorphic to
\[ \Omega^{2i}(M)|_{U}\otimes_{C^{\infty}(U)} \oplus_{i=1}^s  \Gamma(V_{\a'}|_{U_i})\cong 
\Omega^{2i}(M)|_{U}\otimes_{C^{\infty}(U)} \Gamma(V_{\a}|_{U})
\]
which is just $\Omega^{2i}(M,V_{\a})|_U$.
It's clear that these isomorphisms can be patched together and give an isomorphism
$\Omega^{2i}(M,V_{\a})\cong \Omega^{2i}(M',V_{\a'})$ which commutes with $*$ and $d$ since these operators
are defined locally.
$\eta_{\a}(M,s)=\eta_{\a'}(M',s)$, hence $\eta_{\a}(M,0)=\eta_{\a'}(M',0)$

Exactly the same way using the trivial representation for $\a'$ 
one shows that $\eta(M',0)=\eta_{\a(M,M')}(M,0)$.
\end{proof}

In the application we'll have the case that $\pi_1(M')\subset \pi_1(M)$ is normal. We'll now restrict ourselves to this
case. Write $G:=\pi_1(M')/\pi_1(M)$ and write $M_G,\a_G$ for $M'$ and $\a'$. We'll give an explicit matrix representation 
for $\a_G$ and show that if $\a_G\in P_k(\pi_1(M_G))$ then $\a\in P_{ks}(\pi_1(M))$.

Let $g_1,\dots,g_s$ be the elements of $G$ and pick a splitting $\psi:G\to \pi_1(M)$ 
which is of course in general not a homomorphism, but we can arrange $\psi$ such that
$\psi(g^{-1})=\psi(g)^{-1}$ and $\psi(e)=e$.
Let $e_1,\dots,e_k$ denote the
canonical basis of
$\C^k$. Then $g_i \otimes e_j$ is a basis for $\C G\otimes_{\C} \C^k$ and 
$\psi(g_i)\otimes e_j$ is a basis for $\C \pi_1(M)\otimes_{\C \pi_1(M_G)} \C^k$.
We'll write $\a$ with respect to the basis $\psi(g_i)\otimes e_j$.
Note that
\[ a\psi(g_i) \otimes v=\psi(\v(a)g_i)\psi(\v(a)g_i)^{-1}a\psi(g_i)\otimes v=
\psi(\v(a)g_i) \otimes \b(\psi(\v(a)g_i)^{-1}a\psi(g_i)) v \]
since $\v(\psi(\v(a)g_i)^{-1}a\psi(g_i))=1$.
Therefore $\a(a)$ is given by
\[ P_{\v(a)} \hspace{-0.1cm} \bp \a_G(\psi(g_1\v(a))^{-1}a\psi(g_1)) \hspace{-0.1cm}&0&\dots &\hspace{-0.1cm} 0 \\
0\hspace{-0.1cm}&\a_G(\psi(g_2\v(a))^{-1}a\psi(g_2))&&\hspace{-0.1cm}0\\ \vdots\hspace{-0.1cm}
&&\ddots&\hspace{-0.1cm}\vdots\\ 0\hspace{-0.1cm}&0&\dots&\hspace{-0.1cm}\a_G(\psi(g_s\v(a))^{-1}a\psi(g_s))\ep \]
where $P_{\v(a)}:\C G\to \C G$ denotes the matrix corresponding to left multiplication by
$\v(a)$, i.e. $P_{\varphi(a)}(\psi(g_i)\otimes e_j)=\psi(\v(a)g_i)\otimes e_j$.  

Assume that $\a_G$ factors through a homomorphism $\ti{\a}_G$ to a 
  group $H_G$, then $\a$ factors through 
\[ H:=\{ (\v(a),(h\mapsto (\ti{\a}_G(\psi(h\v(a))^{-1}a\psi(h)))_{h\in G}) | a\in \pi_1(M)\} \subset G \ltimes
H_G^G
\]
 where $H_G^G=\mbox{Maps}(G,H_G)$ and $G$ acts on $H_G^G$ by precomposition by left multiplication.   
If $H_G$ is a $p$--group and $G$ is a $p$--group with respect to the same prime, then $H$ is also a $p$--group
which proves the following lemma.

\begin{lemma} \label{lemmaalphaprimetoo}
Let $p$ a prime number. If $G$ is a $p$-group and  $\a_G$ factors through a $p$-group 
then $\a$ also factors through a $p$-group, i.e. $\a \in P_{ks}(\pi_1(M))$.
 \end{lemma}

\section{Computation of eta-invariants for boundary links} \label{sectionetaforboundary}
Let $(L,\varphi)$ be a boundary link pair and $V=V_1 \cup \dots \cup V_m$ a corresponding Seifert surface. 
 Let $\a \in R_k(F_m)$, then define $\t:=\a \circ \v: \pi_1(M_L) \to F_m
\to U(k)$. In this section we'll compute $\rho(M_L,\v)(\a)=\eta_{\t}(M_L)$
using theorem \ref{apsthm}.

First we add handles
$D_i^{2q}
\times D^2$ along the $L_i$ to $D^{2q+2}$  and denote this manifold by $N_L$, then $M_L=\partial(N_L)$.
Note that $\v$ does not extend over $N_L$ since in fact $\pi_1(L)=1$.
We push the surfaces $V_i$ into $D^{2q+2}$, more explicitely,
we can find a map $\iota:V \times I \to D^{2q+2}$, $I=[0,1]$, such that $\iota | V \times 0$ is the embedding of
$V$ into
$S^{2q+1}$, $\iota | L_i \times I$ is constant on the intervals
and such that $\iota|_{int(V)\times I}$ is an embedding.
Now let $\S_i:=\iota(V \times 1) \cup D_i \times 0\subset N_L$, 
$\S:=\cup_{i=1}^m \S_i$, and $N:=N_L
\setminus N(\S)$, then  $\partial(N)=M_L \cup -\S \times S^1$.

We can find embeddings $g_i:D^{2q} \times I \hookrightarrow D_i^{2q} \times D^2$ such 
that $g_i | D^{2q}\times 0$ is just the embedding in $D_i^{2q} \times 0 \subset D_i^{2q}\times D^2$ and such
that  $g_i | D^{2q}\times 1
\subset M_L$ and
$g_i | \partial(D^{2q})\times I \subset V_i$. Now let $T_i:=\iota(V \times I) \cup H_i$ and $T:=\cup_{i=1}^m
T_i$.  Then the Pontrjagin construction for $T \subset N$ gives a map $\pi_1(N) \to F_m$ which extends
$\v:\pi_1(M_L) \to F_m$. We denote the map
$\pi_1(N) \to F_m \to U(k)$ by
$\t$ as well.  Note that $T_i$ inherits an orientation from $\mbox{int}(V_i)\times I \subset T_i$.
By theorem \ref{apsthm}
\[ \eta_{\t}(M_L)-\sum_{i=1}^m \eta_{\ti{\t}_i}(\S_i \times S^1 )=\sign_{\t}(N)- k \cdot \sign(N) \]
where $\ti{\t}_i=\t \circ i_*:\pi_1(\S_i \times S^1) \to \pi_1(N) \to F_m \to U(k)$.

\subsection{Computation of $\eta_{\ti{\t}_i}(\S_i\times S^1)$}
Note that $S^1$ inherits an orientation from the orientations of $\S_i$ and $\S_i \times S^1$.
Denote by $m_i$ the (oriented) generator of $\pi_1(S^1)$,
then
\[   \ti{\t}_i:  \pi_1(\S_i)\times \pi_1(S^1)  \cong \pi_1(\S_i \times S^1)\to U(k)  \]
is given by sending $(g,m_i^e)$ to $\a(t_i)^e$.
We need the following proposition.

\begin{proposition}\cite[thm. 1.2]{N79}
\bn
\item Let $\a_N:\pi_1(N^{2r})\to U(k_N)$ and $\a_X:\pi_1(X^{2s-1})\to U(k_X)$ be representations, then
\[ \eta_{\a_N \otimes \a_X}(N^{2r}\times X^{2s-1})=(-1)^{rs}\sign_{\a_N}(N)\eta_{\a_X}(X) \]
\item Let $\a:\pi_1(S^1)=\Z \to U(1)$ be a representation. If $\a(1)=e^{2\pi ia}, a\in [0,1)$, then
\[ \eta_{\a}(S^1)=\eta(\a(1)):=\left\{ \ba{rl} 0 &\mbox{ if }a=0 \\ 1-2a &\mbox{ if }a\in (0,1) \ea \right. \]
\en
\end{proposition}
 
Therefore
\[ \eta_{\ti{\t}_i}(\S_i \times S^1 )= \eps \sign(\S_i)\sum_{i=1}^m\sum_{i=1}^k\eta(c_{ij}) \]
where  $\{c_{ij}\}_{j=1,\dots,m}$  denotes the set of eigenvalues of $\a(t_i)$ and $\eps:=(-1)^q$.
We can express $\sign(\S_i)$ in terms of the Seifert matrix as follows:
\[ \sign(\S_i)=\sign(V_i)=\sign(\sqrt{\eps}(A_{ii}+\eps A_{ii}^t)) \]
In the case $\eps=-1$ one can easily show that $A_{ii}- A_{ii}^t$ is congruent to 
$\bp 0 & \id \\ -\id &0 \ep$, hence the signature is zero.

\subsection{Computation of $\sign_{\t}(N)$}
\subsubsection{Computation of $H_{q+1}^{\t}(N,\C^k)$}
Denote by $\ti{N}$ the $F_m$-cover of $N$ induced by $\v$, note that
 $C_*(\tilde{N})$ has a right $F_m$-module structure.
Recall that the twisted homology $H_i^{\alpha}(N,\C^k)$ is defined as $H_i(C_*(\tilde{N})\otimes_{\Z F_m}\C^k)$,
where $\C^k$ is a left $F_m$-module via $\a$.  Fix
an orientation preserving embedding
$f: (T,\partial(T))
\times [-1,1]
\to (N,\partial(N))$, such that $f(T \times 0)$ is the usual embedding of $T \subset N$. Let
$X:= N \setminus f(T \times (-1,1))$, then $X$ is homoemorphic to $N$ cut along $T$. We can embed $T $ in $X$ via
the  embeddings $f_+(c) := f(c,1)$ and $f_-(c):=f(c,-1)$.Then $\ti{N} \cong X \times
F_m / \sim $, where $f_-(c_i) \times zt_i \sim f_+(c_i) \times z$ for $c_i \in T_i,z \in F_m$.
From this decomposition we get the following   short exact sequence (where $c_i \in C_*(T_i)$)
\[ \ba{rcccl}
0 \to C_*(T\times F_m)&\to &C_*(X \times F_m)  & \to &C_*(\tilde{N} )
\to 0 \\
(c_i,z) & \mapsto & (f_-(c_i),zt_i)-(f_+(c_i),z) &&\\
&& (c,z) &\mapsto & (c,z) \\
\ea
\]
We tensor with $\C^k$ over $\Z F_m$, the tensored sequence is still exact
since $C_*(\ti{N})$ is a free $\Z F_m$-module.
Taking the long exact homology sequence we get
\[ \dots \to H_i^{\t}(T,\C^k) \to H_i^{\t}(X,\C^k)
\to H_i^{\t}(N,\C^k) \to H_{i-1}^{\t}(T ,\C^k) \to \dots  \]
where
\[ \ba{rcl} H_i^{\t}(T,\C^k)&=&H_i(C_*(T \times F_m) \otimes_{\Z F_m}\C^k)=H_i(C_*(T ) \otimes_{\Z
}\C^k)=H_i(T ,\C^k) \\
    H_i^{\t}(X,\C^k)&=&H_i(C_*(X \times F_m) \otimes_{\Z F_m}\C^k)=
H_i(C_*(X) \otimes_{\Z }\C^k)=H_i(X,\C^k) \ea
\] We have to compute $H_*(X)$. Write $X=X_1 \cup X_2$
where $X_1:=X \cap D^{2q+2}$ and $X_2:=X \cap (\cup_{i=1}^m
D_i^{2q} \times D^2)$.
$X_1$ is homotopy equivalent to a point since
$X_1 \cong D^{2q+2} \setminus f(T\times (-1,1))$, which is just a deformation retract of $D^{2q+2}$.
Furthermore $H_*(X_2)=H_*(\cup_{i=1}^m (D_i^{2q} \times D^2 \setminus H_i))= H_*(m\mbox{ points})$, so from the
Mayer-Vietoris sequence we get for $i \ge 2$
\[  H_i(X) \cong H_{i-1}(X_1 \cap X_2) = H_{i-1}(L \times (D^2 \setminus I))=H_{i-1}(L) \]
furthermore
\[ 0 \to H_1(X) \to H_0(X_1 \cap X_2) \to H_0(m \mbox{ points})\oplus H_0(X_2) \to H_0(X) \to 0 \]
so $H_i(X)=0$ for all $i=1,\dots,2q-1$, $H_0(X)=\Z$ and $H_{2q}(X)=\Z^m$.

\begin{proposition}
If $q>1$ or $(\a(t_i)-\id)$ is invertible for all $i$, then
\[ H_{q+1}^{\t}(N,\C^k) \cong H_{q}(T \times I,\C^k) \cong H_{q}(\S,\C^k) \cong H_{q}(V,\C^k)\]
\end{proposition}

\begin{proof}
The last isomorphism follows since $\S^{2q}=V^{2q} \cup D^{2q}$, the second isomorphism is clear, so 
it only remains to prove the first isomorphism.
For $q \ge 2$ this follows immediately from the long exact sequence. In the case $q=1$ we get 
\[ \dots \to H_2(T,\C^k) \to H_2(X,\C^k)  \to H_2^{\t}(N,\C^k) \to H_{1}(T ,\C^k) \to H_1(X,\C^k) =0  \] 
but the map $H_2(T,\C^k) \to H_2(X,\C^k)$ is induced by the map
\[ \ba{rcl} 
 C_2(T\times F_m) \otimes_{\Z F_m} \C^k &\to &C_2(X \times F_m) \otimes_{\Z F_m} \C^k  \\
(c_i,z)\otimes v & \mapsto & (f_-(c_i),zt_i)\otimes v-(f_+(c_i),z)\otimes v)
\ea
\]
Consider the maps
\[ \ba{rrcccccl} 
 f_+,f_- :& \Z=H_2(\S) & \cong &H_2(T)& \to& H_2(X) &\xrightarrow{\cong} &H_1(X_1 \cap H_2) =\Z \\
&  [ \S ] &\to &[\S] &\to& [f_{\pm}(\S)] &\to & [f_{\pm}(\S) \cap (X_1 \cap X_2)]
\ea \]
But  $[f_{\pm}(\S) \cap (X_1 \cap X_2)]=K$, i.e. $f_+=f_-$.
Therefore
\[ \ba{rcl} 
 C_2(T\times F_m) \otimes_{\Z F_m} \C^k &\to &C_2(X \times F_m) \otimes_{\Z F_m} \C^k  \\
(c_i,z)\otimes v & \mapsto & (f_-(c_i),z)\otimes (t_iv-v)
\ea
\]
If $(\a(t_i)-\id)$ is invertible for all $i$, then $ H_{q+1}^{\t}(N,\C^k) \cong H_{q}(V,\C^k)$.
\end{proof}

In the following we will assume that the assumptions of the proposition hold.  

We can give a more explicit definition of the isomorphism $ H_{q}(V,\C^k)\to H_{q+1}^{\t}(N,\C^k)    $.
Denote by $^*$ the pushing of $V$ into $T=\iota(V\times I)$. Let $\{l_{i1},\dots,l_{ih_i}\}$ be 
 bases of the torsion free part of
$H_{q}(V_i)$, fix  representatives of $l_{ij}$ in $C_{q}(V)$ which are in general position, we'll denote them
by
$l_{ij}$ as well.
For $l \in \{l_{ij},l^*_{ij}\}$ denote by $c^+(l)$ resp. $c^-(l)$ a chain in $X_1 \subset D^{2q+2}$ with
$\partial(c^+(l))=f_+(l)$ respectively $\partial(c^-(l))=f_-(l)$,  we can assume that the 
chains are in general position to each other. 
Consider $X_1$ as lying in $\ti{N}$ via $X_1 \to X_1 \times {e} \to \ti{N} $.
Then the map
\[ \ba{rcl} \psi: C_i(V)\otimes_{\Z F_m} \C^k &\to & C_{i+1}(N)\otimes_{\Z F_m} \C^k \\
        l \otimes v &\mapsto & c^+(l_{ki}) \otimes v-c^-(l_{ki})t_k \otimes v
\ea
\]
induces the above isomorphism 
 $H_{q}(V,\C^k) \cong  H_{q+1}^{\t}(N,\C^k)  $.

\subsubsection{The intersection pairing on $H_{q+1}^{\t}(N,\C^k)$}
Consider the equivariant intersection pairing
\[ \ba{rcl} \langle \:,\: \rangle: C_{q+1}(\tilde{N}) \times C_{q+1}(\tilde{N}) & \to & \Z F_m \\
      (c,\tilde{c}) & \mapsto & \sum_{g \in F_m} ((c \cdot g)\cdot \ti{c} )g^{-1} \ea \]
where $(c \cdot g)\cdot \ti{c}$ is the ordinary intersection number, which is 0 for almost all $g$. Note that
$  \langle cg,\ti{c} \rangle= \langle c,\ti{c} \rangle g$ 
and $ \langle c,\ti{c}g  \rangle=g^{-1} \langle c,\ti{c}  \rangle$. Denote by 
$A$ the Seifert matrix of
$(L,\v)$ with respect to the basis $\{l_{i1},\dots,l_{ih_i}\}$. 

\begin{lemma}
We get the following  matrix of $\langle \:, \: \rangle $ with respect to the elements $\psi(l_{ij})$ 
\[\ba{rcl} &&  \bp A_{11}(1-t_1^{-1})-\eps A_{11}^t(1-t_1 ) & A_{12}(1-t_1)(1-t_2^{-1}) & \dots \\
A_{21}(1-t_2)(1-t_1^{-1}) &A_{22}(1-t_2^{-1})-\eps A_{22}^t(1-t_2)  \\
\vdots & & \ddots \
 \ep =\\
&=& A -\eps    TA^tT^{-1} -AT^{-1}   + \eps  TA^{t}=\\
&=&(A+\eps TA^t)(1-T^{-1})\ea
\]
\end{lemma}

Note that a  similar computation has been done by Ko (cf. \cite{K89}) for the intersection form
of the (abelian) $\Z^m$-cover of $N$.


\begin{proof}
 
Denote by $_+$ resp. $_-$ pushing into
the positive resp. negative direction in $\mbox{int}(V) \times [-1,1] \subset S^{2q+1}$ and in $\mbox{int}(V)
\times I \subset T
\times I$.
 The map 
$\psi:H_{q}(V,\C^k) \to H_{q+1}^{\t}(N,\C^k)$ is induced by
$\psi(l_{ki})=c^+(l_{ki})-c^-(l_{ki})t_k$.  We can deform $\psi(l_{ki})$ into
$d(\psi(l_{ki}))=c^+(l_{ki}^*)-c^-(l_{ki}^*)t_k$.
Note that
\bn
\item  Right multiplication by $t_k$ is an isometry.
\item $\lk(l,\ti{l})=-\eps \lk(\ti{l},l)$.
\item  Denote by $ c^+(l_{ki}),
c^+(l_{lj}^*) \subset X \subset S^{2n+1} \setminus N(L)$ submanifolds representing the corresponding chains.  We
can add a cylinder in $T \times (-1,1)$ to $c^+(l_{lj}^*)$ to get a  submanifold $c$ with $\partial(c)=l_{lj}$.
Then 
$\lk(l_{ki+},l_{ij})=\lk(f_+(l_{ki}),l_{lj})=c^+(l_{ki}) \cdot c= c^+(l_{ki}) \cdot
c^+(l_{lj}^*)$ by the definition of the linking pairing (cf. \cite{R90}).
\item  $ c^-(l_{ki}) \cdot
c^-(l_{lj}^*)=\lk(l_{ki-},l_{lj})$  as above.
\item $c^{\pm}(l) z\cdot c^{\pm}(\ti{l})\ti{z}=0$ if $z \ne \ti{z}$ and if $l,\ti{l}$ don't intersect,  since the
chains  don't intersect.
\item $c^+(l)\cdot c^-(\ti{l})=\lk(l_+,\ti{l}_-)$, $c^-(l)\cdot c^+(\ti{l})=\lk(l_-,\ti{l}_+)$ since  the
embedding $X \subset S^{2n+1} \setminus N(L)$ doesn't change the intersection numbers,  and $f_{\pm}(l)=l_{\pm}$
and 
$f_{\pm}(\ti{l})=\ti{l}_{\pm}$.
\en
Using this we compute
\[ \ba{rcl} \psi(l_{ki})\cdot\psi(l_{lj})&=&\psi(l_{ki})
\cdot d(\psi(l_{lj}))=(c^+(l_{ki})-c^-(l_{ki})t_k) \cdot (c^+(l_{lj}^*)-c^-(l_{lj}^*)t_l)= \\
 &=& c^+(l_{ki}) \cdot c^+(l_{lj}^*)+(-c^-(l_{ki})t_k) \cdot (-c^-(l_{lj}^*)t_l)=\\
&=&\lk(l_{ki+},l_{lj})+\lk(l_{ki-},l_{lj})\delta_{kl}
\ea \]
and for $z \ne 1$ we compute
\[ \psi(l_{ki})z\cdot\vph(l_{lj})=(c^+(l_{ki})z-c^-(l_{ki})zt_k) \cdot (c^+(l_{lj})-c^-(l_{lj})t_l) \]
this is zero except for the following cases:
\[ \ba{rclcl} 
 z&=&t_l& \Rightarrow &\psi(l_{ki})z\cdot\psi(l_{lj})= c^+(l_{ki})t_l \cdot
(-c^-(l_{lj})t_l)=-\lk(l_{ki+},l_{lj}) \\
  z&=&t_k^{-1}& \Rightarrow&\psi(l_{ki})z\cdot\psi(l_{lj})=
 -c^-(l_{ki})t_kt_k^{-1} \cdot  c^+(l_{lj})=-\lk(l_{ki-},l_{lj}) \\
  z&=&t_l^{-1}t_k&  \Rightarrow &\psi(l_{ki})z\cdot\psi(l_{lj})=
 -c^-(l_{ki})t_lt_l^{-1}t_k \cdot (-c^-(l_{lj})t_l)=\lk(l_{ki},l_{lj}) \\
&&&&   \mbox{ since } z \ne 1 \mbox{ implies } k \ne l
\ea \]
The lemma now follows immediately from the definition of the Seifert matrix $A$.
\end{proof}

Recall that the twisted intersection pairing is defined as follows 
\[ \ba{rcl} (,): C_{q+1}(\tilde{N}) \otimes_{\Z F_m}\C^k
\times C_{q+1}(\tilde{N}) \otimes_{\Z F_m}\C^k& \to &\C \\
      (c \otimes v,\tilde{c} \otimes \ti{v}) & \mapsto & \bar{\ti{v}}^t \a( \langle c,\ti{c}  \rangle)  v \ea \]
We now proved the following proposition.
\begin{proposition}
If $(\a(t_i)-\id)$ is invertible for all $i$, then the intersection pairing on $H_{q+1}^{\a}(N,\C^k)$
with respect to the basis $l_{ij} \otimes e_k \in H_{q+1}^{\t}(N,\C^k)$
is represented by the matrix
\[ \sqrt{-\eps}( A -\eps  \a(T) A^t\a(T)^{-1} -A\a(T)^{-1}   + \eps  \a(T)A^{t} \]
In particular 
\[ \sign_{\t}(N)=\sign(\sqrt{-\eps}( A -\eps  \a(T) A^t\a(T)^{-1} -A\a(T)^{-1}   + \eps  \a(T)A^{t}) \]
\end{proposition}

\subsection{Proof of theorem  \ref{thmbdyetacomp}}
Recall that we have to show the following.

\begin{claim} 
Let $(L\subset S^{2q+1},\v)$ be an $F_m$-link, 
$A=(A_{ij})_{i,j=1,\dots,m}$ a Seifert matrix for $(L,\v)$,
 $\a:F_m\to U(k)$ a representation.
Let  $\eps:=(-1)^{q+1}$, then 
\[ \ba{rcccl} \rho(M_L,\v)(\a)&=&
&\eps\sum_{i=1}^m \sign(\sqrt{\eps}(A_{ii}+\eps A_{ii}^t))\sum_{i=1}^m \sum_{j=1}^k \eta(z_{ij})
+\\
&&+&\sign(\sqrt{-\eps}( A -\eps  \a(T) A^t\a(T)^{-1} -A\a(T)^{-1}   + \eps  \a(T)A^{t})) \ea \]
\end{claim}

\begin{proof}
The statement under the assumption  that either $q>1$ or $(\a(t_i)-\id)$ is invertible for all $i$ follows
immediately from the calculations above and the observation that the untwisted signature is 0.

In the general case we have
\[ H_{q+1}^{\t}(N,\C^k) \cong  H_{q}(V,\C^k) \oplus \Im( H_2(X,\C^k)  \to H_2^{\t}(N,\C^k))  \]
Let $ (c,z) \in \Im( H_2(X,\C^k)  \to H_2^{\t}(N,\C^k)) $ for $i=1,2$ and $(d,w) \in 
\Im(H_{q}(V,\C^k) \to H_{q+1}^{\t}(N,\C^k))$.
Then $c g \cdot d=0$ since $c$ can be represented by an element which is supported on $\partial(N)$ whereas
$d$ can be represented by an element which is supported on $N \setminus \partial(N)$.

Since $H_2(X)$ is generated by $[f_-(\S_i)]$ it remains to show 
that $f_-(\S_i) g$ and $f_-(\S_j)^*$ are  disjoint for any $g \in F_m, i,j \in \{1,\dots,m\}$. 
That's obvious for $i \ne j$ and  for $g\ne e$.  Recall that $\S_i \cap (X_1 \cap X_2)=K$.
 Pick a longitude $K'$
for $K$. Pick a Seifert surface $V'$ for $K'$ and close it by a disk $D' $ in the 2-handle over $K$. Then $[V'
\cup D']$ represents
$[\S]$ and we can assume that
$\S_i$ and
$V'
\cup D'$ are in general position. But
\[ \S \cdot (V' \cup D')=(V \cup D^2) \cdot (V' \cup D')=V \cdot V' + D^2 \cdot D'=0\]
since $V\cdot V'=\lk(K,K')=0$ and $D,D'$ can be chosen to be disjoint.

\end{proof}


\begin{thebibliography}{10}
\bibitem[APS75]{APS75}
M. Atiyah, V. Patodi, I. Singer, {\em Spectral asymmetry and Riemannian geometry II}, Math. Proc. Camb. Phil. Soc, 78:
405-432
(1975)
\bibitem[CG86]{CG86}
A. Casson, C.  Gordon, {\em Cobordism of classical knots},
Progr. Math., 62, \'A la recherche de la topologie perdue, 181--199,
 Birkh\"auser Boston, Boston, MA
(1986)
\bibitem[CK99]{CK99}
J. C.ÊCha, K. H. Ko, {\em Signature invariants of links from irregular covers and non-abelian covers},
 Math. Proc. Cambridge Philos. Soc. 127, no. 1: 67--81 (1999)
\bibitem[CO90]{CO90}
T. Cochran, K. Orr, {\em Not all links are concordant to boundary links},
Bull. Amer. Math. Soc. (N.S.) 23, no. 1: 99--106 (1990)
\bibitem[CO93]{CO93}
T. Cochran, K. Orr, {\em Not all links are concordant to boundary links}, Ann. of Math. (2) 138,
no. 3: 519--554 (1993)
\bibitem[CS80]{CS80}
S. Cappell, J. Shaneson, {\em Link cobordism}, Comment. Math. Helv., 55: 20-49 (1980)
\bibitem[F03]{F03}
S. Friedl, Ê{\em Algorithm for finding boundary link Seifert matrices}, preprint (2003)
\bibitem[GL92]{GL92}
P. Gilmer, C. Livingston,
{\em The Casson-Gordon invariant and link concordance},
Topology 31, no. 3: 475--492 (1992)
\bibitem[H67]{H67}
B. Huppert, {\em Endliche Gruppen. I.}, Die Grundlehren der Mathematischen Wissenschaften, Band 134 (1967)
\bibitem[J97]{J97}
D. L. Johnson, {\em Presentations of groups}, Second edition. London Mathematical Society Student Texts, 15.
Cambridge University Press, Cambridge (1997)
\bibitem[K84]{K84}
K. H. Ko, {\em Seifert matrices and boundary links}, Thesis, Brandeis University (1984)
\bibitem[K85]{K85}
K. H. Ko, {\em Algebraic classification of simple links of
odd dimension $\geq$ 3},  Unpublished (1985)
\bibitem[K87]{K87}
K. H. Ko, {\em Seifert matrices and boundary link cobordism},
Transactions of the AMS 299: 657-681 (1987)
\bibitem[K89]{K89}
K. H. Ko, {\em A Seifert-matrix interpretation of Cappell and Shaneson's approach to link cobordisms}, Math.
Proc. Camb. Phil. Soc. 106: 531-545 (1989)
\bibitem[L69]{L69}
J. Levine, {\em Knot cobordism groups in codimension two},
Commentarii Mathematici Helvetici 44: 229-244 (1969)
\bibitem[L69b]{L69b}
J. Levine, {\em Invariants of knot cobordism}, Inventiones Mathematicae 8: 98-110 (1969)
\bibitem[L88]{L88}
J. Levine, {\em Link concordance}, Algebra and topology 1988, 57--76, Korea Inst. Tech., Taeju on (1988)
\bibitem[L89a]{L89a}
J. Levine, {\em  Link concordance and algebraic closure of groups}, Comment. Math. Helv. 64, no. 2: 236--255
(1989)
\bibitem[L89b]{L89b}
J. Levine, {\em  Link concordance and algebraic closure II}, Invent. Math. 96, no. 3: 571--592 (1989)
\bibitem[L90]{L90}
J. Levine, {\em Algebraic closure of groups}, Contemporary Mathematics 109: 99--105 (1990)
\bibitem[L94]{L94}
J. Levine, {\em Link invariants via the eta invariant}, Comment. Math. Helv. 69, no. 1: 82--119 (1994)
\bibitem[L77]{L77}
C.ÊLiang, {\em An algebraic classification of some links of codimension two}, Proc. Amer. Math. Soc. 67, no. 1:
147--151 (1977)
\bibitem[M77]{M77}
T. Matumoto, {\em On the signature invariants of a non-singular complex sesquilinear form}, J. Math.
Soc. Japan 29, no. 1: 67--71 (1977)
\bibitem[M57]{M57}
J. Milnor, {\em Isotopy of links}, Algebraic geometry and topology.
A symposium in honor of S. Lefschetz: 280--306.
Princeton University Press, Princeton, N. J. (1957)
\bibitem[M87]{M87}
W. Mio, {\em On boundary-link cobordism}, Math. Proc. Camb. Phil. Soc., 101: 259-266 (1987)
\bibitem[N79]{N79}
W. Neumann, {\em Signature related invariants of manifolds-I}, Topology 18: 147-172 (1979)
\bibitem[R90]{R90}
D. Rolfsen,  {\em Knots  and links}, Publish or Perish (1990)
\bibitem[S02]{S02}
D. Sheiham, {\em  Invariants of Boundary Link Cobordism}, Preprint
\bibitem[S89]{S89}
L. Smolinsky, {\emÊ Invariants of link cobordism}, Proceedings of the 1987 Georgia Topology Conference
(Athens, GA, 1987). Topology Appl. 32, no. 2: 161--168 (1989)
\bibitem[S89b]{S89b}
L. Smolinsky, {\em A generalization of the Levine-Tristram link invariant},
Trans. Amer. Math. Soc. 315, no. 1: 205--217 (1989)
\bibitem[S65]{S65}
J. Stallings, {\em Homology and central series of groups}, Journal of Algebra 2: 170--181 (1965)
\bibitem[S77]{S77}
N.ÊStoltzfus, {\em Unraveling the integral knot concordance group}, Mem. Amer. Math. Soc. 12, no. 192 (1977)
\bibitem[T73]{T73}
H. F. Trotter, {\em On S-equivalence of Seifert matrices}, Invent. Math. 20: 173-207 (1973)


 \end{thebibliography}
\end{document}